\documentclass[11pt,reqno]{amsart}
\usepackage{amsmath}
\usepackage{amsthm}
\usepackage{amssymb}
\usepackage{enumerate}
\usepackage{amscd}
\usepackage{color}
\usepackage{pb-diagram}
\usepackage{graphicx}

\theoremstyle{plain}

\newtheorem*{prop*}{Proposition}
\theoremstyle{definition}

\theoremstyle{remark}

\numberwithin{equation}{section}

\newcommand{\x}{\noindent}

\newcommand{\vs}{\vspace{0.2cm}}
\newcommand{\vvs}{\vspace{0.4cm}}
\newcommand{\vsss}{\vspace{0.08cm}}

\begin{document}
\title[A Cylinder Theorem for Real K\"ahler Submanifolds]{A Dajczer-Rodriguez Type Cylinder Theorem for Real K\"ahler Submanifolds}\thanks{ Research supported by the NSFC, Grant No. 11271320, and Faculty Stimulus Fund, MRI of OSU}
\author{Jinwen Yan}
\address{Center for Mathematical Sciences\\
                    Zhejiang University\\
                    Hangzhou, 310027 China}
\email{yimkingman@gmail.com}
\author{Fangyang zheng}

\address{Center for Mathematical Sciences\\Zhejiang University \\
Hangzhou, 310027 China, and  Department of Mathematics\\
                    The Ohio State University\\
                    231 West 18th Avenue, Columbus, OH 43210\\}
\email{zheng@math.ohio-state.edu}

\begin{abstract} In 1991, Dajczer and Rodriguez proved in \cite{DR91} that a complete minimal real K\"ahler submanifold of codimension $2$, if with complex dimension $>2$, would be either holomorphic, or a cylinder, or complex ruled. In this article, we generalize their result to real analytic complete  real K\"ahler submanifolds of codimension $4$. The conclusion is that such the submanifold, if with complex dimension $>4$, would be either partially holomorphic, or a cylinder, or a twisted cylinder in the sense that the complex relative nullity foliation is contained in a strictly larger holomorphic foliation, whose leaves are cylinders. We also examine the question of when such a submanifold is complex ruled.
\end{abstract}
\maketitle

\section{Introduction and statement of results}

\vvs

Throughout this article, by a {\em real K\"ahler submanifold} we shall mean  an isometric immersion $f: M^{2n} \rightarrow {\mathbb R}^{2n+p}$ from a K\"ahler manifold of complex dimension $n$ into the Euclidean space. When the codimension is relatively small, the presence of complex structure tends to have strong restrictions on the structure of the submanifold. A typical example of this kind is the classic theorem of Dajczer and Rodriguez in 1991 \cite{DR91}, in which they proved that for any complete minimal $f$ with $p=2$ and $n>2$,  $f$ is either a holomorphic hypersurface, or a cylinder, or complex ruled. To be precise, they proved the following:

\vs

\x {\bf Theorem (Dajczer-Rodriguez \cite{DR91}).} {\em Let $f: M^{2n} \rightarrow {\mathbb R}^{2n+2}$ be a complete minimal real K\"ahler submanifold with $n>2$. Then we have at least one of the following:

 (1) $f$ is holomorphic under an identification ${\mathbb C}^{n+1}\cong {\mathbb R}^{2n+2}$ of the target.

 (2) $f=f_1\times \iota$, where $f_1: N^4 \rightarrow {\mathbb R}^6$ is a complete minimal real K\"ahler submanifold and $\iota$ is the identity map ${\mathbb C}^{n-2}\cong {\mathbb R}^{2n-4}$.

 (3) $M$ admits a holomorphic foliation, whose leaves are totally geodesic flat K\"ahler hypersurface isometric to ${\mathbb C}^{n-1}$, and are mapped by $f$ onto (translations of) linear subspaces in ${\mathbb R}^{2n+2}$.}

\vs

Later in \cite{DG95}, examples and precise descriptions of non-holomorphic minimal real K\"ahler submanifolds in type (3) were given by Dajczer and Gromoll using the Weierstrass representations.

The main purpose of this article is to generalize the above theorem to the codimension $3$ and $4$ cases, and also to allow $f$ to be non-minimal. Before we state our results, we need to recall some well-known facts and also fix some terminologies.

\vsss

Let $f: M^{2n} \rightarrow {\mathbb R}^{2n+p}$ be a real K\"ahler submanifold. Let $\Delta $ be the kernel distribution of the second fundamental form of $f$ and let $\Delta_0=\Delta \cap J\Delta $ be its $J$-invariant part. Denote by $M'\subset M$ the open subset where $\Delta$ reaches its minimum dimension,  $\nu$, which is called the {\em relative nullity index}. In $M'$, $\Delta $ becomes a  foliation with flat totally geodesic leaves. Each leaf will be complete (and isometric to ${\mathbb R}^{\nu }$) when $M$ is assumed to be complete. Similarly, we will denote by $M''\subset M'$ the open subset where $\Delta_0$ reaches its minimum dimension, denoted as $2\nu_0$. $\nu_0$ will be called the {\em complex relative nullity index}. In $M''$, $\Delta_0$ becomes a complex foliation whose leaves are flat, totally geodesic K\"ahler submanifolds in $M$, and when $M$ is complete, each leaf of $\Delta_0$  is holomorphically isometric to $\mathbb{C}^{\nu_0}$.

\vsss

We will call the translations of linear subspaces in ${\mathbb R}^{2n+p}$ {\em linear subvarieties}.

\vs

\x {\bf Definition.} (1). $f$ is {\em partially holomorphic} (or $f$ admits a {\em K\"ahler extension}), if there exists an open dense subset $M_0\subset M$, such that for each connected component $U$ of $M_0$,  $f|_U = h \circ \sigma $ for some real K\"ahler submanifold $h: N^{2n+2k} \rightarrow {\mathbb R}^{2n+p}$ with $k>0$ and some holomorphic embedding $\sigma : U \rightarrow N$.

\vsss

(2). $f$ is called a {\em twisted cylinder}, if there exist an open dense subset $M_0\subset M''$ and a holomorphic foliation ${\mathcal K}$ in $M_0$, such that each leaf of ${\mathcal K}$ consists of parallel leaves of $\Delta_0$, and is of dimension strictly between $2\nu_0$ and $2n$.

\vsss

(3). $f$ is {\em complex ruled,} if there exist an open dense subset $M_0\subset M''$ and a holomorphic foliation ${\mathcal R}$ in $M_0$, such that each leaf of ${\mathcal R}$ consists of parallel leaves of $\Delta_0$, and each leaf is a totally geodesic flat K\"ahler submanifold isometric to ${\mathbb C}^k$ with $\nu_0 < k<n$, mapped by $f$ onto a linear subvariety.

\vs

The terminology in (1) comes from \cite{DG97} and \cite{YZ}, and the one in (2) was coined by Al Vitter in the context of developable submanifolds, see \cite{WZ}. Each leaf of ${\mathcal R}$ in (3) will be called a {\em twisted ruling} for $f$. Note that when $k=n-1$ and $f$ is real analytic, the foliation ${\mathcal R}$ can be extended to the entire $M$ since the limiting position at each boundary point must be unique (see \cite{WZ} for details). In this case $M$ is the total space of the holomorphic vector bundle over a Riemann surface and each fiber is mapped by $f$ onto a linear subvariety.

\vsss

Now we are ready to state our first result:

\vs

\x {\bf Theorem 1.} {\em Let $f: M^{2n} \rightarrow {\mathbb R}^{2n+p}$ be a complete real K\"ahler submanifold with $p\leq 4$ and $f$ being real analytic. We assume that either $n>p$, or $n\leq p$ but $\nu_0>0$. Then either $f$ is partially holomorphic, or a cylinder $f=f_1\times \iota$, where $f_1:N^{2n\!-\!2\nu_0} \rightarrow {\mathbb R}^{2n-2\nu_0+p}$ and $\iota$ the identity map of ${\mathbb C}^{\nu_0} \cong {\mathbb R}^{2\nu_0}$,  or $f$ is a twisted cylinder. }

\vs

Note that when $p=4$, $f$ being partially holomorphic would mean that either $f$ is holomorphic under an identification ${\mathbb C}^{n+2}\cong {\mathbb R}^{2n+4}$, or there exists an open dense subset $M_0\subset M$, such that for each connected component $U$ of $M_0$, there exist a real K\"ahler submanifold $h: N^{2n+2} \rightarrow {\mathbb R}^{2n+4}$ and a holomorphic embedding $\sigma : U\rightarrow N$ so that $f|_U = h\circ \sigma $ holds. Similarly, if $p=3$, then it means the latter case only (while for $p=2$ it means $f$ is holomorphic under an identification ${\mathbb C}^{n+1} \cong {\mathbb R}^{2n+2}$).

\vs

We remark that the real analyticity assumption in the Theorem 1 is a technical one. It is satisfied when $f$ is minimal. Without the real analyticity assumption, the result still holds on an ($\Delta_0$-saturated) open subset of $M$. We added this assumption here mainly to exclude the potential `gluing phenomenon' from happening, namely, cylinders with different ruling directions can potentially be glued together along flat (or flatter) pieces, to form submanifolds which are locally cylinders (with complete leaves) but not global cylinders. It is generally believed that for real K\"ahler submanifolds, such phenomenon does not occur, since the complex structure is too fragile and disallow procedures such as partition of unity. However, at this point, we do not know how to exclude it, even in much simpler situations.

\vsss

As a quick example to illustrate this subtlety, let us recall the conjecture which states that {\em any complete real K\"ahler submanifold with $\nu_0=n-1$ must be a cylinder.} In this case, it is easy to see that each connected component of the open subset $M''$ (where $\Delta_0$ has constant minimum dimension) is a cylinder ${\mathbb C}^{n-1}\times N^2$, but it is still not known (although conjectured so) whether $M$ itself must be a cylinder. (The only known case is when the codimension $p=1$, which implies $\nu_0\geq n-1$. In this case $f$ is a cylinder by the work of Florit and the second author \cite{FZ-hypersurface}). Of course if we assume that $f$ is real analytic, then the local product structure would imply that $f$ is a global cylinder.

\vs

For the sake of convenience, let us introduce the following terminology, in honor of Marcos Dajczer, who pioneered the study of real K\"ahler submanifolds.

\vs

\x {\bf Definition.} A {\em Dajczer submanifold} is a real analytic, complete, real K\"ahler submanifold $f: M^{2n} \rightarrow {\mathbb R}^{2n+p}$ that is neither partially holomorphic nor a cylinder. $f$ is said to be {\em proper Dajczer} if it is Dajczer but not complex ruled.

\vs

With these terminologies at hand, we can restate Theorem 1 simply as: when $p\leq 4$, if  $n>p$ (or $n\leq p$ but $\nu_0>0$), then a Dajczer submanifold is always a twisted cylinder.

\vs

Our next goal is to analyze proper Dajczer submanifolds and see when can they occur.  Let us write $r=n-\nu_0$ and call it the {\em rank} of $f$. When $\Delta_0$ is holomorphic (which will be the case when $f$ is either minimal or complete), its {\em twisting tensors} (see the next section for definition) $C_T$ are all complex linear. Denote by $l$ the maximum (complex) rank of $C_T$ for all $T$ in $\Delta_0$ and we will call $l$ the {\em twisting rank} of $f$. When $f$ is complete, each $C_T$ is nilpotent, so $l$ is an integer between $0$ and $r-1$.

\vsss

Note that $l=0$ when and only when all leaves of $\Delta_0$ are parallel to each other, that is, when $f$ is a cylinder.

\vs

\x {\bf Theorem 2.} {\em Let $f: M^{2n} \rightarrow {\mathbb R}^{2n+p}$ be a real analytic complete real K\"ahler submanifold with $\nu_0>0$. If the twisting rank $l=1$, then $f$ is complex ruled.}

\vs

In particular, if $r=2$, then either $l=0$ and $f$ is a cylinder, or $l=1$ and $f$ is complex ruled. This is proved by Dajczer and Rodriguez (Theorem 2,  \cite{DR91}) when $f$ is minimal and not holomorphic.

\vsss

When the codimension $p=3$ or $4$, the condition $r>p$ would imply that $f$ is partially holomorphic. This is proved by Dajczer-Gromoll \cite{DG97} for $p=3$ and by the authors \cite{YZ} for $p=4$. So we can restate Theorem 2 in the following way:

\vs

\x {\bf Theorem 2*.} {\em Let $f: M^{2n} \rightarrow {\mathbb R}^{2n+p}$ be a real analytic complete real K\"ahler submanifold. Assume either $\nu_0>0$, or $p\leq 4$ and $n>p$. If $f$ is not partially holomorphic, and the twisting rank $l=1$, then $f$ is complex ruled.}

\vs

\vsss

For a complete, real analytic real K\"ahler submanifold $f: M^{2n} \rightarrow {\mathbb R}^{2n+2}$ with $n>2$ that is not minimal, it was proved in \cite{FZ-codim2} that $f$ must be a cylinder. In codimension $3$, we also show that:

\vs

\x {\bf Theorem 3.} {\em Let $f: M^{2n} \rightarrow {\mathbb R}^{2n+3}$ be a real analytic, complete real K\"ahler submanifold which is neither partially holomorphic nor a cylinder. Assume either $n>3$, or $n\leq 3$ but $\nu_0>0$. If $f$ is not minimal, then $f$ is complex ruled.}

\vs

In other words, when $n>3$ (or $\nu_0>0$), any proper Dajczer submanifold in codimension $3$ must be minimal.  It would be rather interesting to give precise description of this type of submanifolds, in the spirit of \cite{DG95}. See also \cite{Hennes} for related discussions.

\vsss

For $p=4$, the non-minimal case also tend to be more restrictive than the minimal case. For instance, we will show that $l=3$ can only occur when $f$ is minimal. Most of the $l=2$ subcases (in fact all but one) will lead to $f$ being complex ruled. So proper Dajczer submanifolds in codimension $4$ again form a rather special and restrictive class. See Theorem 4 and related discussion in the next section for more details.

\vvs

\noindent {\bf Acknowledgement:} We are grateful to Marcos Dajczer and his collaborators for their inspiring papers on the subject of real K\"ahler submanifolds, which opened the way to the studies of this area in submanifold theory. The second author would like to thank his former collaborators Luis Florit and Wing San Hui. The present work is a continuation of these earlier joint works as well as \cite{YZ}. We would also like to thank CMS of Zhejiang University for the support, and in particular to Hongwei Xu for his help and interest. The second author is also supported by a Faculty Stimulus Grant from MRI of Ohio State University and the National Science Foundation of China.

\vvs

\vvs

\vvs

\section{The proof of the theorems}

\vvs

Throughout this section, we will assume that $f: M^{2n}\rightarrow {\mathbb R}^{2n+p}$ is a  real analytic, complete, real K\"ahler submanifold. Denote by $M''\subset M'\subset M$ and $\Delta$, $\Delta_0$, $\nu$, $\nu_0$ as before. Assume $\nu_0>0$, then both $M''$ and $M'$ are open dense in $M$, and $\Delta_0$ is a holomorphic foliation in $M''$ with flat totally geodesic leaves isometric to ${\mathbb C}^{\nu_0}$, which are mapped by $f$ onto linear subvarieties.

\vsss

For $x\in M''$, denote by $\Delta_0^{\perp }$ the orthogonal complement of $\Delta_0$ in the tangent space $T_xM$. Recall that for any $T\in \Delta_0$, the {\em twisting tensor} $C_T: \Delta_0^{\perp } \rightarrow \Delta_0^{\perp }$ is defined by
$$ C_T(Y) = - (\nabla_Y\tilde T)^{\perp } $$
for any $Y\in \Delta_0^{\perp}$. Here $Z^{\perp}$ denote the $\Delta_0^{\perp}$-component of $Z$, and $\tilde T$ is any vector field in a neighborhood of $x$ with value $T$ at $x$.  $C_T$ is well-defined, namely, independent of the choice of $\tilde T$, and is a tensor. It satisfies the equations
\begin{eqnarray}
& & \nabla_S C_T \ = \ C_T C_S + C_{\nabla_S T} \\
& & (\nabla_X C_T)Y - (\nabla_Y C_T)X \ = \ C_{(\nabla_XT)^{\Delta_0} }Y - C_{(\nabla_YT)^{\Delta_0}} X
\end{eqnarray}
for any $T,S$ in $\Delta_0$ and any $X$, $Y$ in $\Delta_0^{\perp }$. We refer the readers to \cite{DG90} and \cite{DR91} for the history and more details on this tensor. When $f$ is minimal, it is proved in \cite{DG} that $C_T$ is complex linear, namely, $JC_T = C_T J$ where $J$ is the almost complex structure of $M$. When $f$ is assumed to be complete, it was proved in \cite{FZ-codim2} (see also \cite{WZ-JDG} and \cite{WZ-slight}) that $C_T$ is always complex linear. When $f$ is complete, it is also well-known that all $C_T$ are nilpotent.

\vs

Following \cite{WZ}, define distributions ${\mathcal K}$ and ${\mathcal R}$ in $M''$ respectively by
$$ {\mathcal K}_x = \Delta_0+ \bigcap \mbox{ker}(C_T) , \ \  {\mathcal R}_x = \Delta_0+ \sum \mbox{Im}(C_T), \ \ {\mathcal R}^0_x={\mathcal K}_x\cap {\mathcal R}_x $$
where the intersection and sum are taking over all $T$ in $\Delta_0$ and $\mbox{Im}(C_T)$ denotes the image space of $C_T$. In an open dense subset $M_0\subset M''$, ${\mathcal K}$, ${\mathcal R}$, and ${\mathcal R}^0$ reach their minimum dimensions thus become subbundles of the tangent bundle of $M$. All three are $J$-invariant since all $C_T$ are complex linear. We claim that:

\vs

\x {\bf Lemma 1.} {\em ${\mathcal K}$ is a holomorphic foliation in $M_0$. }

\vs

\x {\em Proof.}  First let us show that ${\mathcal K}$ is a foliation. Given any $X$, $Y$ in $\Delta_0^{\perp}$ with $C_TX=C_TY=0$ for all $T$, and any $S\in \Delta_0$, we want to show that $[S,X]$ and $[X,Y]$ are both in ${\mathcal K}$. For convenience, we can extend $C_T$ to $TM$ by letting it act trivially on $\Delta_0$ and still denote it as $C_T$. Apply (2.1) to $X$, we get $C_T(\nabla_SX) =0 $. So $C_T([S,X])=0$. Similarly, apply (2.2) to $X$ and $Y$ would imply $C_T(\nabla_XY - \nabla_YX)=C_T([X,Y])=0$. So ${\mathcal K}$ is a foliation.

To see that ${\mathcal K}$ is holomorphic, we need to show that for any tangent vector $X$ and any $Y$ in ${\mathcal K}$, the vector $\nabla_XY+\nabla_{JX}JY$ is always in ${\mathcal K}$. It suffices to consider the case when $Y\in \Delta_0^{\perp}$, since $\Delta_0$ is known to be holomorphic. Apply (2.2) to $X$, $Y$, and then also to $JX$, $Y$, we get
\begin{eqnarray*}
C_T(\nabla_X Y) & = & \nabla_Y (C_T X) - C_T (\nabla_Y X) + C_{ (\nabla_YT)^{\Delta_0} } X \\
 C_T(\nabla_{JX} Y) & = & \nabla_Y (C_T JX) - C_T (\nabla_Y JX) + C_{ (\nabla_YT)^{\Delta_0} } JX
\end{eqnarray*}
Since both $\nabla$ and $C_T$ commute with $J$, these two equations lead to the fact that $C_T(\nabla_XY + \nabla_{JX}JY)=0$. \qed

\vs

It is easy to see that ${\mathcal K}$ is parallel along each leaf of $\Delta_0$, and the leaves of ${\mathcal K}$ consists of parallel leaves of $\Delta_0$ thus are $\Delta_0$-cylinders.

\vs

Next we consider the distribution ${\mathcal R}$ and ${\mathcal R}^0$. In general ${\mathcal R}$ may not be a foliation, but ${\mathcal R}^0$ is always a holomorphic foliation. Note that we have ${\mathcal R}^0_x=\Delta_x^0\oplus {\mathcal L}_x$, where ${\mathcal L}$ consists of all $X\in \Delta_0^{\perp }$ such that $X=\sum C_{S_i}Z_i$ is a finite sum for some $S_i\in \Delta_0$ and  $Z_i \in \Delta_0^{\perp }$, and $C_TX=0$ for all $T$ in $\Delta_0$.

\vs

\x {\bf Lemma 2.} {\em ${\mathcal R}^0$ is a holomorphic foliation.  }

\vs

\x {\em Proof.} First let us show that ${\mathcal R}^0$ is a foliation. Let $X$, $Y$ be in ${\mathcal L}$ and $S$, $T$ in $\Delta_0$. We want to show that $[S,X]$ and $[X,Y]$ both lie in ${\mathcal R}^0=\Delta_0\oplus {\mathcal L}$. We already know by Lemma 1 that ${\mathcal K}$ is a foliation, so it suffices to show that they lie in ${\mathcal R}$. For $[S,X]$, since $\nabla_XS \in \Delta_0$ as $C_SX=0$, it suffices to show $\nabla_SX \in {\mathcal R}$. Write $X=\sum C_{S_i}Z_i$, and apply (2.1) to $Z_i$, we get
$$ \nabla_SC_{S_i}Z_i \equiv 0 \ \  \mbox{mod} \ {\mathcal R}, $$
so $\nabla_SX \in {\mathcal R}$, thus $[S,X]\in {\mathcal R}^0$. For $[X,Y]$, apply (2.2) to $Y$ and $Z_i$, and use the fact that $C_TY=0$ for any $T$, we get
$$ \nabla_YC_{S_i}Z_i \equiv 0 \ \ \mbox{mod} \ {\mathcal R},$$
so $\nabla_YX\in {\mathcal R}$. Similarly, $\nabla_XY\in {\mathcal R}$, thus $[X,Y]\in {\mathcal R}$, so ${\mathcal R}^0$ is a foliation.

Next we show that ${\mathcal R}^0$ is holomorphic. It suffices to prove that for any $X$ and any $Y\in {\mathcal L}$, $\nabla_XY + \nabla_{JX}JY $ lies in ${\mathcal R}$. Write $Y=\sum C_{S_i}Z_i$, and apply (2.2) to $X$, $Z_i$ and $JX$, $Z_i$ respectively, we get
\begin{eqnarray*}
\nabla_X C_{S_i}Z_i & \equiv & \nabla_{Z_i} (C_{S_i} X) \ \ \mbox{mod} \ {\mathcal R} \\
\nabla_{JX} C_{S_i}Z_i & \equiv & \nabla_{Z_i} (C_{S_i} JX) \ \ \mbox{mod} \ {\mathcal R}
\end{eqnarray*}
so $\nabla_XY+\nabla_{JX}JY \equiv 0$ mod ${\mathcal R}$, and ${\mathcal R}^0$ is holomorphic. \qed

\vs

In the particular case when ${\mathcal R}\subset {\mathcal K}$, the holomorphic foliation ${\mathcal R}^0={\mathcal R}$ is actually a complex ruling:

\vs

\x {\bf Lemma 3.} {\em Assume that ${\mathcal R}\subset {\mathcal K}$. Then the holomorphic foliation ${\mathcal R}^0={\mathcal R}$ is actually a complex ruling, namely, each leaf of ${\mathcal R}$ is a totaly geodesic K\"ahler submanifold of $M^n$, holomorphically isometric to ${\mathbb C}^k$ for some $k\geq v_0$, and is mapped by $f$ onto a linear subvariety.}

\vs

\x {\em Proof.} If suffices to show that $\widetilde{\nabla}_XY \in {\mathcal R}$ for any vector fields $X$ and $Y$ in ${\mathcal R}$, where $\widetilde{\nabla}$ is the connection in ${\mathbb R}^{2n+p}$. Notice that the assumption ${\mathcal R}\subset {\mathcal K}$ means that $C_TC_S=0$ for any $T$, $S$ in $\Delta_0$. Write $X=\sum C_{T_i}Z_i$ and $Y=\sum C_{S_k}W_k$. Note that the twisting tensor $C_T$ satisfies the following well-known property:
\begin{eqnarray}
\alpha (C_TZ,W) = \alpha (Z, C_TW)
\end{eqnarray}
for any tangent vectors $Z$ and $W$, where $\alpha$ is the second fundamental form of $f$. From (2.3), we get
$$ \alpha (X, Y) = \sum \alpha (Z_i, C_{T_i}C_{S_k}W_k) = 0,$$
so $\widetilde{\nabla}_XY =\nabla_XY$. Since $C_{S_k}X=0$, by applying (2.2) to $X$ and $W_k$, we get
$$ \nabla_XY = \sum \nabla_X(C_{S_k}W_k) \equiv 0 \ \ \mbox{mod} \ {\mathcal R}, $$
so ${\mathcal R}$ is indeed totally geodesic in ${\mathbb R}^{2n+p}$, thus an open subset of a linear subvariety $P$. By the real analyticity of $f$, we know that the entire $P$ must be contained in $M$. \qed

\vs

Of course the main difficulty is to know when the common kernel of $C_T$ will be non-trivial. To this end we observe that the main algebraic result in \cite{WZ}, Proposition 2, exactly dealt with this question. In fact when $f$ is minimal, the situation can be exactly carried over to here. When $f$ is not minimal, we need to modify things a little bit and prove a similar result.

\vsss

Let us again denote by $H$ and $S$ respectively the $(1,1)$ and $(2,0)$ component of the second fundamental form $\alpha$ of $f$. We will fix a generic point $x\in M$ and examine the algebraic relations. Denote by $V\cong {\mathbb C}^r$ the complex vector space (of type $(1,0)$ tangent vectors) corresponding to $\Delta_0^{\perp }$, and restrict $H$ and $S$ on $V$. Then $H: V\times \overline{V} \rightarrow W$ and $S: V\times V \rightarrow W$ are respectively  Hermitian symmetric  or symmetric bilinear form, where $W\cong {\mathbb C}^p$ is the complexification of $T_xM^{\perp }\cong {\mathbb R}^p$. We will extend the inner product $\langle , \rangle $ complex linearly to $W$. They satisfy the following relations:
\begin{eqnarray}
& & \langle H_{X\overline{Y}}, H_{Z\overline{W}}\rangle \ = \ \langle H_{Z\overline{Y}}, H_{X\overline{W}}\rangle  \\
& & \langle H_{X\overline{Y}}, \ S_{ZW}\rangle \ = \  \langle H_{Z\overline{Y}}, S_{XW}\rangle  \\
& & \langle S_{XY}, \ S_{ZW}\rangle \ = \  \langle S_{ZY}, S_{XW}\rangle \\
& & \ \mbox{ker}(H) \cap \mbox{ker}(S) =0
\end{eqnarray}
for any $X$, $Y$, $Z$, and $W\in V$. Here $\mbox{ker}(H)$ consists of all $X\in V$ such that $ H_{X\overline{Y}} =0$ for all  $Y \in V$, and $\mbox{ker}(S)$ likewise.

\vsss

Now let $A^1, \ldots , A^k$ be complex linear transformations on $V$ such that for any complex linear combination $A=\sum_{i=1}^k t_iA^i$, it holds
\begin{eqnarray}
&& H_{\!AX, \overline{Y}}=0, \ \ \ S_{\!AX,Y}=S_{\!AY, X}, \ \ \ \forall \ X,Y \in V\\
&& A \  \mbox{is} \ \mbox{nilpotent.}
\end{eqnarray}

Denote by $N(A)$ the kernel of $A$, and $R(A)$ the range (image space) of $A$. We claim that:

\vs

\x {\bf Lemma 4.} {\em Let $H$, $S$ and $A$ be as above satisfying (2.4) through (2.9). If $r\leq 4$, then the common kernel space  \
$ \bigcap_{i=1}^k N(A^i) \neq 0$.}

\vs

\x {\em Proof.} If $H=0$, then this is exactly the situation of Proposition 2 in \cite{WZ}, and the proof can be carried over without any change.

In general, let us consider $V'=\mbox{ker}(H)\subset V$. By (2.8), we have $R(A)\subset V'$. Denote by $S'$ the restriction of $S$ on $V'\times V'$. If $X\in \mbox{ker}(S')$, then for any $Y\in V$, we have
$$ S_{AX, Y}=S_{X,AY}=0$$
since $AY \in V'$. So $AX$ is contained in $\mbox{ker}(S)$. By (2.7), we get $AX=0$ for any $A$ so $X$ is in the common kernel of $A^i$.

Therefore, we may assume that $\mbox{ker}(S')=0$. Now consider $A'=A|_{V'}$. Since $R(A)\subset V'$, we have $A': V'\rightarrow V'$. Clearly, it is nilpotent and is symmetric with respect to $S'$, so we are once again in the $H=0$ situation. Note that $N(A')=N(A) \cap V'$, so if $X$ is in the common kernel of $A'$, it  would be in the common kernel of $A$ as well.  \qed

\vs

Now we are ready to prove Theorem 1.

\vs

{\em Proof of Theorem1:} Let $f$ be as in Theorem 1. Since $p\leq 4$, by the main theorem of \cite{YZ} (the $p=3$ case is due to Dajczer and Gromoll in \cite{DG97}), we know that $f$ must be partially holomorphic unless $r\leq p$. So when $n>p$, we have $\nu_0=n-r>0$. The other assumption is $n\leq 4$ and $\nu_0>0$. Thus,  if we assume that $f$ is not partially holomorphic, then under the dimension assumption of Theorem 1, we would always have $\nu_0>0$. So we have a holomorphic, totally geodesic foliation $\Delta_0$, with $r\leq 4$. By Lemma 4, we know that the common kernel of all $C_T$ is non-trivial, so by Lemma 1 we get a strictly larger holomorphic foliation ${\mathcal K}$ in some open dense subset of $M$. The leaves of ${\mathcal K}$ are $\Delta_0$-cylinders. If the dimension of the leaves of ${\mathcal K}$ are $2n$, then $f$ is a cylinder. Otherwise it is a twisted cylinder, and this completes the proof of Theorem 1. \qed

\vs

Next let us recall that the {\em twisting rank} $l$ of $f$ is defined to be the maximum of the complex rank of $C_T$ for all $T$ in $\Delta_0$. Since $f$ is complete, all $C_T$ are nilpotent, so $0\leq l\leq r-1$, where $r=n-\nu_0$ is the rank of $f$. Although $f$ might not be a holomorphic immersion, the completeness of $f$ and the holomorphicity of $\Delta_0$ works just like in the situation of developable holomorphic submanifolds of ${\mathbb C}^N$ in \cite{WZ}, and we have the following:

\vs

\x {\bf Lemma 5.} {\em If the twisting rank $l=1$, then ${\mathcal R} \subset {\mathcal K}$.}

\vs

\x {\em Proof.} When $f$ is minimal, we are exactly in the situation of Case 1 in the proof of Proposition 2 in \cite{WZ}, and the proof there can be carried over without any change. For the general case, since $l=1$, each $A^i$ is of rank $1$ (we skip the trivial ones), so each kernel $N(A^i)$ is a codimension one subspace in $V$, and each range space $R(A^i)$ is one dimensional. Since each $A^i$ is nilpotent, we have $R(A^i) \subset N(A^i)$ for each $i$. So if all $N(A^i)$ coincide, then $A^iA^j=0$ for any $i$, $j$. Now let us assume that not all the kernel spaces coincide, say, $N(A^1)\neq N(A^2)$. Choose a basis $\{ e_1, \ldots , e_r\}$ of $V$, such that $\{ e_1, \ldots , e_{r-1}\}$ is a basis of $N(A^1)$ and $\{ e_1, \ldots , e_{r-2}, e_r\}$ is a basis of $N(A^2)$. We claim that $R(A^1)=R(A^2)$ in this case. Note that $R(A^1)$ is spanned by $A^1(e_r)$, while $R(A^2)$ is spanned by $A^2(e_{r-1})$. Consider the matrix  $A=A^1+A^2$. We have $A(e_{r-1})=A^2(e_{r-1})\in R(A^2)$ and $A(e_r)=A^1(e_r)\in R(A^1)$. Since $A$ has rank at most one, we know that $R(A^1)=R(A^2)$ must hold.

For any $i>2$, $N(A^i)$ must be different with either $N(A^1)$ or $N(A^2)$. Thus $R(A^i)$ must be equal to either $R(A^1)$ or $R(A^2)$, so the range space of all $A^i$ are the same, which implies $A^iA^j=0$ for all $i,j$. This completes the proof of Lemma 5. \qed

\vs

{\em Proof of Theorem 2 and 2*:} Note that Theorem 2 is the immediate consequence of the combination of Lemma 3 and Lemma 5. For Theorem 2*, as we noticed before, we always have $\nu_0>0$  under the dimension assumptions, so Theorem 2 applies. \qed

\vs

 {\em Proof of Theorem 3:} Let $f: M^{2n} \rightarrow {\mathbb R}^{2n+3}$ be a complete, real analytic, real K\"ahler submanifold that is neither partially holomorphic nor a cylinder. That is, we have a Dajczer submanifold of codimension $3$. The dimension assumption is either $n>3$, or $n\leq 3$ but $\nu_0>0$. In the former case we have $\nu_0>0$ by \cite{DG97}. So we have the complex relative nullity foliation $\Delta_0$ whose leaves have complex dimension $r\leq 3$. When $r=1$ or $2$, we already know that $f$ must be a cylinder or complex ruled, thus we may assume that $r=3$. The twisting rank $l$ is either $1$ or $2$. If $l=1$, $f$ is again complex ruled by Theorem 2. So we just need to show that the case $l=2$ cannot occur unless $f$ is minimal.

 In fact, when $l=2$ and $f$ is not minimal, we will show that $f$ must be partial holomorphic, which is ruled out by the assumption.  We will use notations in the discussion right before Lemma 4. Note that we are now in the special case that $r=p=3$ and $H\neq 0$. Since $l=2$, we may take a $A$ with rank $2$. By (2.8), $R(A) \subset \mbox{ker}(H)$, so $R(A)=\mbox{ker}(H)=V' \cong {\mathbb C}^2$. Let $A'=A|_{V'}: V'\rightarrow V'$. It is also nilpotent, and non-trivial. So there exist a basis $\{e_2, e_3\}$ of $V'$ such that $A(e_2)=e_3$ and $A(e_3)=0$. Extend it into a basis $\{ e_1, e_2, e_3\}$ of $V$, and write $A(e_1)=ae_2+be_3$. Then $a\neq 0$ since $A$ has rank $2$. Replace $e_1$ by $\frac{1}{a}e_1$, we may assume that $a=1$.

 By (2.8), $A$ is symmetric with respect to $S$, so we have $S_{33}=S_{23}=0$ and $S_{13}=S_{22}$. Here we wrote $S_{ij}$ for $S_{e_ie_j}$. Notice that $S_{13}\neq 0$, since otherwise $e_3$ would be in both $\mbox{ker}(S)$ and $\mbox{ker}(H)$ at the same time.

 For $H$, the only non-zero entry would be $H_{1\overline{1}}$, which we will denote by $\lambda_1\xi^1$ where $\xi^1$ is a unit vector and $\lambda_1>0$. By (2.5) and (2.6), we know $\langle \xi^1, S_{13}\rangle =0$, and $\langle S_{13}, S_{13}\rangle =0$ since $S_{33}=0$. So $S_{13}=\lambda(\xi^2+\sqrt{-1}\xi^3)$ where $\lambda >0$ and $\{ \xi^1, \xi^2, \xi^3\}$ forms an orthonormal frame of $TM^{\perp } \cong {\mathbb R}^3$. The symmetry relation (2.5) and (2.6) imply that $\langle \xi^2 , S\rangle = -\sqrt{-1} \langle \xi^3, S\rangle $. So as in \cite{YZ}, we see that the subbundle $E\subset TM^{\perp }$ spanned by $\{ \xi^2, \xi^3\}$ admits an almost complex structure $J$, which will imply that $f$ is partially holomorphic.  This completes the proof of  Theorem 3. \qed

 \vvs

 Now we turn our attention to the codimension $4$ case. Let $f$ be a Dajczer submanifold in codimension $4$, namely, $f: M^{2n}\rightarrow {\mathbb R}^{2n+4}$  is a real analytic complete real K\"ahler submanifold which is neither partially holomorphic nor a cylinder. We assume that either $n>4$, or $n\leq 4$ but $\nu_0>0$. Then by \cite{YZ} we have $\nu_0>0$ and $r\leq 4$. Theorem 1 says that $f$ is always a twisted cylinder. But we would like to know when will $f$ be complex ruled, namely, when will the complex relative nullity foliation $\Delta_0$ be contained in a bigger holomorphic foliation whose leaves are linear subvarieties. Or equivalently, we want to know when can $f$ be a proper Dajczer submanifold.

 \vsss

 If $r=2$, then by Theorem 2 of \cite{DR91}, we know that $f$ is complex ruled, regardless of the codimension. When $r=3$, Theorem 3 states that, if $f$ has codimension $3$, then $l$ must be $1$ (thus $f$ is complex ruled) unless $f$ is minimal. That proof does not work when the codimension is $4$, as $S$ can have a fourth component which is of rank $2$ to mess up the partial holomorphicity. However, if $f$ is non-minimal and has $l=2$, then the normal bundle of $f$ admits a rather special decomposition, which might lead to structural results on such submanifolds. When $r=4$, we can show that the case $l=3$ can only occur when $f$ is minimal, while $l=2$ would also imply that $f$ is `close' to being minimal. To make it precise, we need the following terminology:

 \vs

 \x {\bf Definition.} For a real K\"ahler submanifold $f: M^{2n}\rightarrow {\mathbb R}^{2n+p}$, its {\em $H$-index}, or {\em non-minimality index}, is defined to be the maximum (amongst all points in $M$) of the real dimension of the image space of $H$, the $(1,1)$-component of the second fundamental form of $f$. We will denote it as $h_f$ or simply $h$.

 \vs

 When $f$ is complete, real analytic, and with $\nu_0>0$, by Corollary 9 of \cite{FZ-codim2} we know that $f$ will be a cylinder if $h\geq p-1$, where $p$ is the codimension. Following the line of argument as in the proof of Theorem 3, we have:

 \vs

 \x {\bf Theorem 4.} {\em Let $f: M^{2n}\rightarrow {\mathbb R}^{2n+4}$ be a real analytic complete real K\"ahler submanifold with $n>4$ (or $n\leq 4$ and $\nu_0>0$) such that $f$ is neither partially holomorphic nor a cylinder. If $l=3$ then $f$ must be minimal. If its $H$-index $h\geq 2$, then $l=1$ and $f$ is complex ruled. }

 \vs

 Note that under the assumption of Theorem 4, we have $\nu_0>0$, $r\leq 4$, and $f$ is always a twisted cylinder. Theorem 4 says that $f$ will be complex ruled unless it is minimal ($h=0$) or `almost minimal' (in the sense that the $h=1$). We omit the proof here since it is strictly analogous to that of Theorem 3. It would be very interesting to further analyze the structure of those $f$ that are not complex ruled, namely, the proper Dajczer submanifolds in codimension $p\leq 4$.

\vvs

\vvs

\vvs

\vvs

\end{document}